%
%
%  AUTHORS:  FRANC FORSTNERIC
%
%  TITLE:    OKA'S PRINCIPLE FOR MULTIVALUED SECTIONS OF RAMIFIED MAPPINGS
%
%  FINAL VERSION:  November 12, 2001. 
%  
%_____________________________________________________________________________
%
%   This is the Macro file.
%
%
\scrollmode
\magnification=\magstep1
\parskip=\smallskipamount
\hoffset=1.5cm \hsize=12cm

\def\demo#1:{\par\medskip\noindent\it{#1}. \rm}
\def\ni{\noindent}               % noindent
     % new indented par in horizontal mode
 % new unindented par in horizontal mode
\def\ll{\leftline}
\def\cl{\centerline}

\def\begin{\ll{}\vskip 10mm\nopagenumbers}  % beginning of the paper
\def\pn{\footline={\hss\tenrm\folio\hss}}   % pagenumbers at bottom
\def\ii#1{\itemitem{#1}}

%
%   \beginsection
%
\outer\def\beginsection#1\par{\bigskip
  \message{#1}\leftline{\bf\&#1}
  \nobreak\smallskip\vskip-\parskip\noindent}

%
%  \proclaim
%
\outer\def\proclaim#1:#2\par{\medbreak\vskip-\parskip
    \noindent{\bf#1.\enspace}{\sl#2}
  \ifdim\lastskip<\medskipamount \removelastskip\penalty55\medskip\fi}

\def\endpr{\hfill $\spadesuit$ \medskip}

%
%
%  Roman capital characters
%
%
               % roman D in math mode

               % roman T
               % roman N in math

%
%  Special letters: real, complex numbers etc.
%

\def\R{{\rm I\kern-0.2em R\kern0.2em \kern-0.2em}}
\def\N{{\rm I\kern-0.2em N\kern0.2em \kern-0.2em}}
\def\P{{\rm I\kern-0.2em P\kern0.2em \kern-0.2em}}
\def\B{{\rm I\kern-0.2em B\kern0.2em \kern-0.2em}}
\def\C{{\rm C\kern-.4em {\vrule height1.4ex width.08em depth-.04ex}\;}}
\def\CP{\C\P}

%
%
% Caligraphic capital characters
%
%

\def\cH{{\cal H}}
\def\cI{{\cal I}}
\def\cJ{{\cal J}}

\def\cL{{\cal L}}

\def\cO{{\cal O}}

\def\cS{{\cal S}}
\def\cT{{\cal T}}

\def\cVT{{\cal VT}}

%
%
%  Small Greek letters in Math mode
%

\def\d{\delta}
\def\e{\epsilon}

\def\l{\lambda}

\def\c{\chi}
%\def\o{\omega}

%
%
%  Miscellaneous symbols
%
%
\def\bar{\overline}              % conjugate
\def\bs{\backslash}              % backslash
\def\di{\partial}                % partial derivative
\def\dibar{\bar\partial}         % di-bar derivative

%
%  Symbols
%
             % disc
        % closed disc
               % boundary of the disc

\def\c*{{\C}^*}
                         % roman Aut in math mode

\def\wt{\widetilde}
\def\wh{\widehat}

%
% Abbreviations
%
\def\dim{{\rm dim}}                    % dimension
\def\holo{holomorphic}                   % holomorphic
                  % automorphism
                  % homomorphism
               % analytic subset
                % homeomorphism
                    % continuous
\def\nbd{neighborhood}                   % neighborhood
                   % pseudoconvex
        % strongly psc
                   % real-analytic
               % plurisubharmonic

                    % totally real
             % polynomially convex
         % holo convex
           % holomorphic function
\def\ss{\subset\!\subset}                % relatively compact subset
                       % such that
                  % support
              % C^n equivalent

\def\hvf{holomorphic vector field}
\def\hvb{holomorphic vector bundle}

\def\deg{{\rm deg}}

\def\br{{\rm br}}

%_____________THE TEXT____________________________________________

\nopagenumbers
\ll{}
\vskip 1cm
\cl{\bf THE OKA PRINCIPLE FOR MULTIVALUED SECTIONS}
\cl{\bf OF RAMIFIED MAPPINGS}
\bigskip
\cl{\bf Franc Forstneri\v c} 
\bigskip\medskip\rm

%
%
%  INTRODUCTION
%
%
\beginsection 0. Introduction.

A central problem in the analysis of \holo\ mappings $h\colon Z\to X$ between 
complex manifolds (or complex spaces) is to construct \holo\ sections, i.e., 
maps $f\colon X\to Z$ satisfying $h(f(x))=x$ for all $x\in X$. 
If $h$ is unramified (a holomorphic submersion), there exist plenty of local 
\holo\ sections passing through any point of $Z$ and the main question is the 
existence of global sections. The most one may hope for is that every continuous 
sections of $h$ can be homotopically deformed to a \holo\ section. This holds 
only for special classes of submersions and its validity is commonly referred 
to as the {\bf Oka principle} for sections of $h$. The strongest result in this 
direction was given in [F3], extending the works of Oka [O], Grauert [Gr1, Gr2], 
Forster and Ramspott [FR], and Gromov [Gro] (for proofs of Gromov's results 
see [FP1--FP3] and [L]). 

In this paper we consider the analogous problem for {\bf ramified mappings}.
A point $z\in Z$ is a {\bf ramification point} of $h\colon Z\to X$
if $dh_z\colon T_z Z\to T_{h(z)} X$ is not surjective, and the set of
all such points is the {\bf ramification locus} $\br_h$. 
(For maps of complex spaces we include $Z_{\rm sing}\cup h^{-1}(X_{\rm sing})$
into $\br_h$.) There need not 
exist any local sections, not even continuous ones, passing through a 
ramification point (consider $h(z)=z^2$ at $z=0$). Furthermore, the existence 
of a continuous section need not imply the existence of a local \holo\ 
section at such points. For instance, the function $h(z,w)=z^pw^q$ on $\C^2$ 
for coprime $p,q\in \N$ admits the H\"older-continuous section 
$f(re^{i\theta})=\bigl( r^{1/2p}e^{ ip'\theta}, r^{1/2q}e^{iq'\theta}\bigr)$,
where $p',q'$ are integers with $pp'+qq'=1$, but it admits no local 
single-valued \holo\ section with $f(0)=(0,0)$.

The natural objects associated to ramified maps $h\colon Z\to X$ onto $X$
are complex subvarieties $V\subset Z$ such that $h|_V\colon V\to X$ is
a proper finite projection onto $X$, i.e., an {\bf analytic cover} [W].
If the base $X$ is a connected complex manifold or an irreducible complex space, 
there are an integer $d\in\N$ and a closed complex subvariety $\d\subset X$ 
such that for all $x\in X\bs \d$ the fiber $V_x=V\cap h^{-1}(x)$ consists 
of precisely $d$ distinct points while the exceptional fibers $V_x$ over 
points $x\in \d$ contain less than $d$ points ($d$ is the {\bf degree} and 
$\d$ is the {\bf discriminant locus} of the analytic cover $h\colon V\to X$).

If we consider the points in $V_x$ with appropriate algebraic multiplicities
then every fiber $V_x$ consists of precisely $d$ points and  we may 
think of $V$ as the {\bf graph of a \holo\ $d$-valued section} $F$ of $h$.
This notion can be made precise by considering $F$ as a fiber preserving 
\holo\ map from $X$ to the $d$-fold symmetric power $Z^d_{\rm sym}$
of $Z$. The latter point of view is especially useful since we may 
consider the topological objects in $Z$, with combinatorial properties 
analogous to those of analytic covers, as graphs of continuous (or smooth) 
multivalued sections of $h$. The graph of a general holomorphic
multivalued section of $h$ is an {\bf analytic chain} in $Z$ with a finite
proper $h$-projection onto $X$, that is, a formal combination 
$V=\sum_j m_j V_j$ of finitely many closed irreducible complex 
subvarieties $V_j\subset Z$ with coefficients $m_j\in\N$ such that 
$h|_{V_j} \colon V_j\to X$ is an analytic cover for each $j$. 

The main result of the paper, Theorem 1.1, is a version of the Oka principle 
for multivalued sections of ramified maps onto a Stein base which are
elliptic submersions over the complement of the ramification locus.
It applies in particular to maps $h$ whose regular fibers are compact 
complex curves of genus zero or one (Corollary 1.2). We also prove the Oka 
principle for liftings of  \holo\ maps (Theorem 1.3 and Corollary 1.4). 
We hope that these results will be useful in the study of ramified mappings 
which arise naturally in analytic and algebraic geometry.

\pn

\beginsection 1. The results.

All complex spaces  are assumed to be reduced and finite 
dimensional. $X_{\rm sing}$ denotes the singular locus of $X$ and 
$X_{\rm reg}=X\bs X_{\rm sing}$. Let $h\colon Z\to X$ be a \holo\ map
of complex spaces. We denote by $\br_h \subset Z$ its {\bf ramification locus}, 
consisting of all points $z\in Z$ such that $Z$ is singular at $z$ or $X$ 
is singular at $x=h(z)$ or $dh_z\colon T_z Z\to T_x X$ is not surjective.

For $d\in \N$ we denote by $Z^d_{\rm sym}$ the $d$-fold {\bf symmetric power} 
of $Z$, the quotient of the Cartesian power $Z^d$ by the action of the 
permutation group on $d$ elements permuting the entries of 
$(z_1,\ldots,z_d)\in Z^d$. Thus a point in $Z^d_{\rm sym}$ 
is an unordered $d$-tuple of points in $Z$. $Z^d_{\rm sym}$ 
inherits from $Z^d$ the structure of a complex space and hence we 
may speak of continuous resp.\ \holo\ maps $X\to Z^d_{\rm sym}$. 
For details we refer to section 4 below and  to [W].

%
%  d-valued sections
%
\proclaim  Definition 1:  A \holo\ (resp.\ continuous) {\bf $d$-valued section} 
of $h$ is \holo\ (resp.\ continuous) map $F\colon X\to Z^d_{\rm sym}$ such that 
$F(x)\subset h^{-1}(x)$ for each $x\in X$. 
$F$ is {\bf unramified} at $x\in X$ if its restriction to some small open 
\nbd\ $U\subset X$ of $x$ is a union of $d$ continuous (resp.\ holomorphic) 
single-sections of $h$. We denote by $\br_F$ the {\bf ramification locus} of $F$, 
consisting of all points $x\in X$ such that $F$ is not a union of single-valued
sections in any \nbd\ of $x$.

The number $d$ is called the {\bf degree} 
of $F$. Let $\#F(x)$ denote the number of distinct points in 
$F(x)$ and set   
$$ \mu_F= \max\{\#F(x)\colon x\in X\} \le d, \quad
    \d_F=\{x\in X\colon \#F(x)<\mu_F\}.
$$
The set $\d_F$ is called the {\bf discriminant locus} of $F$.
Clearly $\br_F\subset \d_F$ and both sets are closed when $F$
is continuous. If $F$ is \holo\ then $\br_F$ and $\d_F$ are
closed complex subvarieties of $X$.

%
%
%  DEFINITION OF A SPRAY
%
%
\medskip\ni \bf Definition 2. \sl  
A \holo\ map $h\colon Z \to X$ between complex spaces is an 
{\bf elliptic submersion} over an open subset $\Omega\subset X_{\rm reg}$ if 
the restriction $h\colon h^{-1}(\Omega)\to \Omega$ is a submersion 
of complex manifolds and each point $x\in \Omega$ has an open \nbd\ 
$U$ such that there exist a \hvb\ $p\colon E\to Z|_U:=h^{-1}(U)$ 
and a \holo\ map $s\colon E\to Z|_U$ satisfying the following conditions 
for each $z\in Z|_U$: 
\item{(i)} $s(E_z) \subset Z_{h(z)}$ (equivalently, $h s=h p$),
\item{(ii)} $s(0_z)=z$, and 
\item{(iii)}  the derivative $ds \colon T_{0_z} E \to T_z Z$ maps the
subspace $E_z \subset T_{0_z} E$ surjectively onto the vertical tangent 
space $VT_z(Z):=\ker dh_z$.
\rm \medskip

A tripple $(E,p,s)$ as above is called a (fiber-) {\bf dominating spray} 
associated to the submersion $h\colon Z|_U \to U$ (condition (iii) is 
the domination property of $s$). To our knowledge sprays were first 
introduced by M.\ Gromov in [Gro, 1.1.B] as a replacement for the
exponential map on Lie groups which was used in the proof of the classical
Oka-Grauert principle. 

The notation $\cH^k(A)=0$ means that $A$ has the $k$-dimensional 
Hausdorff measure zero; this makes sense for subsets of analytic spaces.

\proclaim 1.1 Theorem:
{\bf (Oka's principle for multivalued sections)}
Let $h\colon Z\to X$ be a \holo\ map of a complex space $Z$ onto an irreducible 
$n$-dimen\-sional Stein space $X$. Assume that $X_0$ is a closed complex subvariety 
of $X$ containing $h(\br_h)$ such that $h$ is an elliptic submersion over 
$X\bs X_0$. Let $F$ be a continuous $d$-valued section of $h$ which is \holo\ in 
a \nbd\ of $X_0$, unramified over $X\bs X_0$
and satisfies $\cH^{2n-1}(\d_F)=0$. Then there is a homotopy 
$F_t\colon X\to Z_{\rm sym}^d$ $(t\in [0,1])$ of continuous $d$-valued sections 
of $h$ such that $F_0=F$, each $F_t$ is \holo\ in a \nbd\ of $X_0$, 
unramified over $X\bs X_0$ and satisfies $F_t(x)=F(x)$ for $x\in X_0$, 
and $F_1$ is \holo\ on $X$.

We have already mentioned that the graph of a \holo\ multivalued section of $h$ 
is a pure $n$-dimensional analytic chain in $Z$ with proper finite projection 
onto $X$ (Proposition 4.2), and hence Theorem 1.1 may be viewed as an existence result 
for such chains. The problem of extending a \holo\ $d$-valued section of $h$ defined
locally near $X_0$ to a continuous $d$-valued section over $X$ can be 
treated by methods of {\it obstruction theory} 
(see for instance [Wd]).

\demo Remarks concerning Theorem 1.1:  
1.\ The space $Z$ need not be Stein (the fibers of $h$ may even 
be compact). Theorem 1.1 is new even for $d=1$ since the known results only 
apply to single-valued sections of unramified maps. Under suitable hypotheses 
the result also holds over a reducible Stein space $X$ (for single-valued sections 
no extra hypotheses are needed).  Theorem 1.1  holds under the weaker 
condition that $h$ is a {\bf subelliptic submersion} over $\wt X=X\bs X_0$ 
(this is explained in the subsequent paper [F3]). The same applies to 
Theorems 1.3 and 2.1 below. 

\ni 2. Theorem 1.1 holds with the usual additions described 
(in the case of  single-valued sections) in [Gro, FP2, FP3]. 
For instance, if $F$ is \holo\ in a \nbd\ of $K\cup X_0$ for 
some compact, holomorphically convex subset $K\subset X$ and if $h$ is elliptic 
over $X\bs (K\cup X_0)$ then the homotopy $F_t$ in Theorem 1.1 may be chosen 
such that each $F_t$ is  \holo\ in a \nbd\ of $K\cup X_0$, it approximates $F$ 
uniformly on $K$, and it agrees with $F$ to a given finite order along $X_0$.  
(See Theorem 1.4 in [FP3] for a precise statement of such a result 
for single-valued sections.) 
\endpr

\ni 3. The condition $\cH^{2n-1}(\d_F)=0$ is satisfied for any
\holo\ multivalued section since $\d_F$ is a proper complex subvariety of $X$.
This condition implies that $\d_F$ is nowhere dense in $X$ and its complement 
$X\bs \d_F$ is path connected and locally path connected (provided that $X$ 
is irreducible) which guaranties a unique decomposition of $F$ into irreducible 
components (Proposition 4.1). This is no longer the case if $\cH^{2n-1}(\d_F)>0$. 
For instance, the 2-valued map $F\colon \C\to \C^2_{\rm sym}$ given by 
$F(x+iy)=[|x|,-|x|]$ (with $\d_F=\{x=0\}$ of Hausdorff dimension one) 
has two splittings into single-valued continuous maps: 
(a) $F_1(z)=|x|$, $F_2(z)=-|x|$, and (b) $F_1(z)=x$, $F_2(z)=-x$. 

\ni 4. A jet-transversality argument shows that for a generic smooth 
perturbation of any multivalued section the set $\d\cap X_{\rm reg}$ 
is a smooth real submanifold of real codimension at least two and hence 
$\cH^{2n-1}(\d)=0$. However, such a generic perturbation may 
introduce additional ramification points. 

\demo Example 1: Let $\chi\colon\C\to\R_+$ be a smooth 
function which vanishes precisely on $D=\{|z| \le 1\}$. 
Consider the maps $F_\e \colon\C\to \C^2_{\rm sym}$ 
defined by
$$ F_\e(z)= [(\chi(z)+\e)\sqrt z, -(\chi(z)+\e) \sqrt z \,], \qquad z\in\C.
$$
Clearly $F_0$ is unramified and $\d_{F_0}=D$, but for any $\e >0$ the map 
$F_\e$ is ramified at $z=0$ and satisfies $\d_{F_\e}=0$. There exists 
no unramified perturbation $G$ of $F_0$ satisfying $\cH^{1}(\d_G)=0$
which can be seen as follows. The normalized graph of $G$ (obtained by separating 
the self-intersections as in the proof of Lemma 5.1) is a covering space over 
$\C$ and hence trivial. This is a contradiction since the graph of $F$ over 
any circle $\{|z|=r\}$ for $r>1$ is a nontrivial covering space 
which remains nontrivial after a small perturbation.
\endpr

\ni\bf 1.2 Corollary. \sl
Let $h\colon Z\to X$ be a \holo\ map of a complex space $Z$ onto an 
irreducible $n$-dimen\-sional Stein space $X$. Assume that $X_0$ is 
a closed complex subvariety of $X$ % containing $X_{\rm sing}$ 
such that $h\colon \wt Z=Z\bs h^{-1}(X_0) \to \wt X=X\bs X_0$ is a submersion 
of complex manifolds.  Then Theorem 1.1 applies in each of the following cases:

\item{(a)} each connected component of the fiber $Z_x$ for $x\in \wt X$ 
is either a rational curve $(\CP^1)$ or an elliptic curve (a complex torus);

\item{(b)} the restriction $h\colon \wt Z\to \wt X$
is locally trivial (a \holo\ fiber bundle) and the fiber $Z_x$ 
is a complex Lie group or a complex homogeneous space;

\item{(c)} $\wt Z= V\bs \Sigma$ where $h\colon V\to \wt X$ 
is a \hvb\ over $\wt X$ of rank $k \ge 2$ and $\Sigma$ is 
complex subvariety of the associated bundle $\wt V\to \wt X$ with fibers 
$\wt V_x\simeq \CP^k$ such that $\dim \Sigma_x\le k-2$ for all $x\in \wt X$.

\demo Proof: In each case the restricted submersion
$h\colon \wt Z\to \wt X$ is elliptic and hence Theorem 1.1 applies.
Case (b) was  considered in [Gr1, Gr2] and [FP2], and case (c) 
was considered in [Gro] and, more explicitly, in [FP2, Theorem 1.7]. 
In case (c) the fibers $\Sigma_x$ are algebraic subvarieties of 
$\wt V_x\simeq \CP^k$. 

In case (a) the connected components of the 
fiber $Z_x=h^{-1}(x)$ for $x\in \wt X$ are all of the same type 
(either $\CP^1$ or elliptic). In the first case the complex 
structure on $Z_x$ is independent of $x$ and hence $\wt Z\to \wt X$ 
is a fiber bundle with complex homogeneous fiber $\CP^1$, so the result 
is a special case of (b). If the components of $Z_x$ are elliptic curves 
$C_{x,j}$ ($1\le j\le j_0$), the parameter of the complex structure on 
$C_{x,j}$ is locally a holomorphic function of $x$. Hence the universal 
covering maps $\C\to C_{x,j}$ can be chosen to be locally holomorphic 
in $x$, and these maps give sprays on $Z|_U$ over small open sets 
$U\subset \wt X$. 
\endpr

\demo Example 2: 
The following is an explicit example of a fibration of type (a) in Corollary 1.2. 
Let $S$ be a compact complex surface in $\CP^N$. Choose a point 
$p\in \CP^N\bs S$ and let $X=\CP^{N-1}$ denote the set of all complex 
hyperplanes $\l\subset \CP^N$ passing through $p$. Set 
$Z=\{ (\l,z)\colon \l\in X,\ z\in \l\cap S \}$
and denote by $h\colon Z\to X$ the projection $h(\l,z)=\l$. The ramification 
locus $\br_h$ is the set of points $(\l,z)\in Z$ such that the intersection 
of $\l$ with $S$ is non-transverse at $z$. Since $h$ is proper, its 
projection $X_0=h(\br_h) \subset X$ is a closed complex subvariety of $X$. 
For $\l\in X \bs X_0$ the fiber $h^{-1}(\l)= \l\cap S$ is a union of finitely
compact Riemann surfaces whose genus $g$ is independent of $\l\in X$. 
If $g=0$ then $S$ is called a {\bf ruled surface}, and if $g=1$ then $S$ is an 
{\bf elliptic surface} (see [BV]). In each of these two cases Theorem 1.1 
holds over any Stein domain in $X$, but it fails when $g\ge 2$ because of
the hyperbolicity.
\endpr

%
% LIFTINGS
%
Our next result extends Theorem 2.1 in [F2] to ramified maps $h$.

\proclaim 1.3 Theorem: {\bf (The Oka principle for liftings)}
Let $h\colon Z\to X$ and $f\colon Y\to X$ be \holo\ maps of complex spaces.
Assume that $X_0\subset X$ is a closed complex subvariety containing
$f(Y_{\rm sing}) \cup h(\br h)$ and $g_0\colon Y \to Z$ is a 
continuous map which is \holo\ in an open set containing $Y_0:= f^{-1}(X_0)$ 
and satisfies $h g_0=f$. If $Y$ is Stein and  $h$ is an elliptic submersion 
over an open \nbd\ of the set $f(Y\bs Y_0)$ in $X$ then for each $k\in\N$
there exists a homotopy of continuous maps $g_t\colon Y\to Z$ such that 
for each $t \in [0,1]$ we have $hg_t=f$, $g_t$ and $g_0$ are tangent to 
order $k$ along $Y_0$, and the map $g_1$ is \holo\ on $Y$.
If in addition $g_0$ is \holo\ in a \nbd\ of a compact holomorphically convex 
subset $K\ss Y$, the homotopy $g_t$ can be chosen such that, in addition to 
the above, it approximates $g_0$ uniformly on $K$.

The following diagram illustrates Theorem 1.3: 
$$ \matrix{ & & Z & \cr
            & {g \atop\nearrow} & \downarrow & \!\!\!\!\!\!\!\! \! h \cr
            Y & {f\atop \longrightarrow} & X & \cr}
$$
A map $g$ for which this diagram commutes is called a {\it lifting} of $f$, 
and the result is that (under the stated conditions) a continuous lifting
can be homotopically deformed to a \holo\ lifting. The spaces $Z$ and $X$ in 
Theorem 1.3 need not be Stein (only $Y$ is  Stein). Theorem 1.3 
is proved in sect.\ 3.

We can apply Theorem 1.3 to the construction of entire maps on vector bundles 
whose images avoid certain complex subvarieties. Let $h\colon E\to X$ be a \hvb\ 
of rank $q$ over a Stein manifold $X$. For each $x\in X$ we denote by 
$\wh E_x \cong \CP^q$ the compactification of the fiber $E_x \cong \C^q$ obtained 
by adding the hyperplane at infinity $\Lambda_x\cong\CP^{q-1}$. The resulting 
fiber bundle $\wh h\colon \wh E \to X$ is again \holo\ since the transition maps
for $E$, which are $\C$-linear automorphisms of fibers $E_x$, 
extend to projective linear automorphisms of $\wh E_x$.

\proclaim 1.4 Corollary:
Let $h\colon E\to X$ be a \hvb\ of rank $q$ over a Stein manifold
$X$ and let $\wh E\to X$ be the associated bundle with fiber $\CP^q$
as above. Let $\Sigma$ be a closed complex subvariety of $\wh E$ whose 
fiber $\Sigma_x$ has complex codimension at least two in $\wh E_x$ and satisfies
$0_x\notin \Sigma_x$ for each $x\in X$. Then for every $k\in\N$ 
there exists a fiber-preserving \holo\ map $F\colon E\to E\bs \Sigma$ 
which is tangent to the identity to order $k$ along the zero section of $E$.

The conclusion of Corollary 1.4 can also be stated as follows:
{\it There exists a family of entire mappings $F_x$ on 
the fibers $E_x\cong \C^q$, depending holomorphically on the parameter $x\in X$, 
such that $F_x$ is tangent to the identity at the origin $0_x\in E_x$ and 
its image $F_x(E_x)$ misses the subvariety $\Sigma_x \subset E_x$ for each $x\in X$.} 

Of course the point $0_x$ can be replaced by $g(x)$ where $g\colon X\to E\bs \Sigma$
is a \holo\ section. Note that each $\Sigma_x$ is projective-algebraic 
by Chow's theorem. The result is false in general if $\Sigma_x$ has 
codimension one in $E_x$ (since its complement may be
Kobayashi hyperbolic).

\demo Proof of Corollary 1.4: Let $Z=E\bs \Sigma$. The hypothesis
on $\Sigma$ implies that the restricted submersion $h|_Z\colon Z\to X$ is 
elliptic (see Corollary 1.8 in [FP2]). Take $Y=E$ (which is Stein), 
let $Y_0$ denote the zero section of $E$ and set $Z=E\bs \Sigma$.
By hypothesis we have $Y_0\cap \Sigma=\emptyset$. Choose 
a continuous fiber preserving map $g_0\colon E\to Z$ which  
equals the identity in an open set $U\subset E$ containing $Y_0$
(such $g_0$ can be obtained by contracting each fiber $E_x$ to 
a \nbd\ of $0_x$ which does not intersect $\Sigma_x$). 
Since $h|_Z\circ g_0=h$, $g_0$ is a continuous lifting of the 
\holo\ map $h \colon Y=E\to X$ with respect to the submersion
$h|_Z \colon  Z\to X$. By Theorem 1.3 we obtain a homotopy 
of liftings $g_t\colon Y\to Z$ from $g_0$ to a \holo\ lifting $g_1\colon Y\to Z$ 
such that the homotopy is fixed to order $k$ along $Y_0$. The map $F=g_1$ 
satisfies Corollary 1.4.
\endpr

\ni\it Open problem: \rm  Is it possible to choose $F \colon E\to E\bs \Sigma$ 
in Corollary 1.4 to be {\it injective} (i.e., such that $F_x$ is a 
Fatou-Bieberbach map on $E_x\simeq \C^q$ for each $x\in X$) ? 
By Proposition 1.4 in [F1] the answer is affirmative 
when $X$ consists of a single point.
\medskip

In the remainder of this section we explain the organization
of the paper. In section 2 we prove Theorem 1.1 for single-valued sections
of ramified maps (this is stated separately as Theorem 2.1). 
The main point is that the proof of Oka's principle for 
elliptic submersions, given in [FP2] and [FP3], extends to ramified
maps provided that one can construct a local spray 
around the graph of any \holo\ section of $h$ over a holomorphically 
convex subset of $X$, such that this local spray is dominating 
outside of the branch locus of $h$. In [FP2, FP3] we used the 
fact that the vertical tangent space $VT(Z)=\ker dh$ of a holomorphic 
submersion $h\colon Z\to X$ is a \hvb\ over $Z$ and hence is generated 
over any open Stein subset $\Omega\subset Z$ by finitely 
many sections (vertical \hvf s). The composition of local flows of these 
sections is a local spray on $\Omega$. 
When $h$ has ramification points, $VT(Z)$ is no longer a 
vector bundle but merely a linear space over $Z$. Nevertheless, germs of 
\holo\ sections of this space form a coherent analytic sheaf over $Z$ 
which is locally free over $Z\bs \br_h$ (Proposition 2.2), and this
enables us to complete the proof. 

In section 3 we reduce Theorem 1.3 to Theorem 2.1. We associate to the
map $f\colon Y\to X$ the pull-back $\wh h\colon \wh Z\to Y$ of 
$h\colon Z\to X$ so that sections of $\wh h$ over $Y$ are in one-to-one 
correspondence with maps $g\colon Y\to Z$ satisfying $hg=f$ (i.e., with 
liftings of $f$). Furthermore, if $U\subset X$ is an open subset of $X$
and $V=f^{-1}(U)\subset Y$ then any $h$-spray associated to $h^{-1}(U)\to U$ 
pulls back to an $\wt h$-spray associated to $\wh h^{-1}(V)\to V$.
Hence Theorem 1.3 follows from Theorem 2.1 applied to sections of $\wt h$.

In section 4 we recall the basic results on multivalued sections
which are used in the proof of Theorem 1.1 for $d>1$. 

In section 5 we deduce the general case of Theorem 1.1 from Theorem 1.3
as follows. Given a $d$-valued section $F\colon X\to Z^d_{\rm sym}$
of $h\colon Z\to X$ as in Theorem 1.1, we construct a normal
Stein space $Y$ and a continuous map $g_0\colon Y \to Z$ 
satisfying the following:
\item{--} the composition $f=hg_0\colon Y\to X$ is a $d$-sheeted
analytic cover onto $X$ which is unramified over $X\bs X_0$, 
\item{--} $g_0$ maps the fiber $Y_x:=f^{-1}(x)$ onto $F(x)$ for 
each $x\in X$, and 
\item{--}  $g_0$ is \holo\ in $f^{-1}(U_0)$ if $F$ is \holo\ 
in $U_0\subset X$.

\ni The space $Y$ should be thought of as the {\bf normalized graph} of 
$F$ in which the multiple points over $X\bs X_0$ have been separated. 
Since $f\colon Y\to X$ is a proper finite map of $Y$ 
onto a Stein space $X$, $Y$ is also Stein. The inverse of $f$ is 
a \holo\ $d$-valued section $f^{-1}\colon X\to Y^d_{\rm sym}$ of $f$. 
The map $g_0\colon Y\to Z$ is a continuous lifting of $f=hg_0\colon Y\to X$ 
with respect $h\colon Z\to X$. Theorem 1.3 provides a homotopy of liftings
$g_t\colon Y\to Z$ of $f$, connecting $g_0$ to a holomorphic lifting $g_1$. 
Then $F_t=g_t f^{-1}\colon X\to Z^d_{\rm sym}$ is a homotopy of $d$-valued 
sections of $h\colon Z\to X$ satisfying Theorem 1.1.

\beginsection 2.\ Proof of Theorem 1.1 for single-valued sections.

In this section we prove the following version of Theorem 1.1 for
single-valued sections. All complex spaces are assumed to be
reduced and finite dimensional.

\proclaim 2.1 Theorem:
Let $h\colon Z\to X$ and $X_0\subset X$ be as in Theorem 1.1
(hence $h$ is an elliptic submersion over $X\bs X_0$).
For any continuous section $F\colon X\to Z$ which is \holo\ in an open 
set containing $X_0$ and for any $k\in\N$ there exists a homotopy 
$F_t \colon X\to Z$ $(t\in[0,1])$ of continuous sections such that 
$F_0=F$, for each $t\in [0,1]$ the section $F_t$ is \holo\ in a \nbd\ 
of $X_0$ and tangent to $F_0$ to order $k$ along $X_0$, 
and $F_1$ is \holo\ on $X$. If $F$ is \holo\ in a \nbd\ of
$K\cup X_0$ for some compact, holomorphically convex subset $K$
of $X$ then we can choose $F_t$ to be \holo\ in a \nbd\ of
$K\cup X_0$ and to approximate $F=F_0$ uniformly on $K$.

When $Z$ and $X$ are complex manifolds and $h$
is a surjective submersion (i.e., $\br_h=\emptyset$), Theorem 2.1 is
a special case of Theorem 1.4 in [FP3]; when $X_0=\emptyset$ it
is included [Gro, 4.5 Main Theorem] and in [FP2, Theorem 1.5].
The presence of ramification points of $h$ requires a refinement
of the proof related to the construction of local sprays in \nbd s of
graphs of \holo\ sections over Stein sets in $X$ which we shall
now describe.

\medskip\ni\bf 2.2 Proposition. (Existence of local sprays)
\sl Let $h\colon Z\to X$ be a \holo\ map of reduced complex spaces.
For any open Stein subset $\Omega$ of $Z$ there exist an integer $N\in \N$,
an open set $V \subset \Omega\times\C^N$ containing $\Omega \times \{0\}^N$,
and a \holo\ map $s\colon V \to Z$ satisfying the following:
\item{(i)}   $s(z,0)=z$ for $z\in \Omega$,
\item{(ii)}  $h(s(z,t))=h(z)$ for $(z,t)\in V$,
\item{(iii)} $s(z,t)=z$ when $(z,t)\in V$ and $z\in \br_h$,
\item{(iv)}  for each $z\in \Omega\bs \br_h$ the derivative at
$t=0\in\C^N$ of $t\to s(z,t)\in Z$  maps $T_0\C^N$
surjectively onto $VT_z Z :=\ker dh_z$.
\medskip\rm

A map $s$ satisfying Proposition 2.2 is called a
{\bf local spray} for $h$ over $\Omega$. 
Indeed $s$ satisfies all properties of a spray except
that it is not defined globally on $\Omega\times \C^N$ and 
the domination property (iv) only holds in the complement of the 
ramification locus $\br_h$. When $h$ is a submersion of complex manifolds, 
Proposition 2.2 coincides with Lemma 5.3 in [FP1].

\demo Proof of Proposition 2.2: Our reference are Chapters
1 and 2 in [Fi]. Recall that the {\bf  tangent space}
of a complex space $Z$ is a linear space $\pi\colon TZ\to Z$ 
(a complex space with linear fibers over $Z$) obtained as follows. 
Fix a point $z_0\in Z$ and represent an open \nbd\
$Z_0\subset Z$ of $z_0$ as a closed complex subspace of an open
subset $W$ of $\C^m$, defined by a sheaf of ideals
$\cI\subset \cO_W$ which is generated by \holo\ functions
$f_1,\ldots,f_r\in \cO(W)$. Thus the structure sheaf
$\cO_{Z_0}$ is isomorphic to the quotient $\cO_W/\cI$
restricted to $\{f=0\}=Z_0$.
If $(w_1,\ldots,w_m,\xi_1,\ldots,\xi_m)$ are coordinates in $W\times \C^m$
then $TZ_0=TZ|_{Z_0}$ is the closed complex subspace
of $W\times \C^m$ generated by the functions
$$ f_1,\ldots,f_r\ \ {\rm and}\ \
   {\di f_i\over \di w_1}\,\xi_1 + \cdots + {\di f_i\over \di w_m}\,\xi_m
   \ \ {\rm for}\ \ i=1,\ldots,r.                        \eqno(2.1)
$$
The projection $TZ_0\to Z_0$ is induced by $W\times \C^m\to W$,
$(w,\xi)\to w$. Different local representations of $Z$ in $\C^m$ give
isomorphic representations  of the tangent space. Over $Z_{\rm reg}$ the
space $TZ$ is the usual tangent bundle of $Z$.

Suppose furthermore that $h\colon Z\to X$ is a \holo\ map of complex
spaces. The {\bf  vertical tangent space} $\pi\colon VT(Z)\to X$ with
respect to $h$ (also called in [Fi] `the tangent space of $Z$ over $X$'
and denoted $T(Z/X)$) is a linear space over $Z$ (a subspace of $TZ$)
with the following local description.
Fix a point $z_0\in Z$ and let $x_0=h(z_0)\in X$.
Let $Z_0\subset W\subset\C^m$ be a local representation of a \nbd\
of $z_0$ as above and let $X_0\subset W'\subset \C^n$ be a local
representation of a \nbd\ $X_0\subset X$ of $x_0$. We may choose
these neighborhoods such that $h(Z_0)\subset X_0$ and the restriction
$h\colon Z_0\to X_0$ extends to a \holo\ map
$H=(H_1,\ldots,H_n)\colon W\to\C^n$. One takes $T(Z/X)|_{Z_0}$
to be the closed complex subspace of $W\times \C^m$ generated by the
functions (2.1) together with
$$ {\di H_i\over \di w_1}\, \xi_1 + \cdots + {\di H_i\over \di w_m}\,\xi_m
   \ \ {\rm for}\ \ i=1,\ldots,n.                      \eqno(2.2)
$$
Again the result is independent of the local representations
up to isomorphism and the complex space obtained in this way coincides
with the usual vertical tangent bundle over $Z\bs \br_h$.
The spaces $TZ$ and $VT(Z)$ need not be reduced even if $Z$ is,
but this will not be important for our purposes.

We denote by $\cT_Z$ (resp.\ $\cVT_Z$) the sheaf of germs of \holo\ sections
of $TZ$ (resp.\ of $VT(Z)$). These  are $\cO_Z$-analytic sheaves
which are free over $Z\bs \br_h$. Sections of $\cT_Z$ are called
{\bf vector fields} on $Z$ and sections of $\cVT_Z$ are called
{\bf vertical vector fields}.

%
%  Coherence of sheaves of vector fields.
%
\proclaim 2.3 Lemma: The $\cO_Z$-analytic sheaves
$\cT_Z$ and $\cVT_Z$ are coherent.

\ni \it Proof. \rm
Indeed the sheaf $\cL$ of holomorphic sections of any linear
space $L\to Z$ over a complex space $Z$ is a coherent
$\cO_Z$-analytic sheaf [Fi, p.\ 53, Corollary]. In the present
case this can be seen directly as follows. Using a local representation
for $VT(Z)|_{Z_0}$ as above, the sheaf $\cVT_{Z_0}$ consists of
germs of \holo\ maps $\xi=(\xi_1,\ldots,\xi_m)\colon Z_0\to \C^m$
satisfying
$$
    \sum_{k=1}^m {\di f_i(w)\over \di w_k}\, \xi_k(w) =0
    \ \ (1\le i\le r);\ \
    \sum_{k=1}^m {\di H_j(w)\over \di w_k}\, \xi_k(w) =0
    \ \  (1\le j\le n).
$$
To get $\cT_Z$ we only take the first set of equations. Thus  both sheaves
are locally sheaves of relations and hence coherent [GR, p.\ 131].
\endpr

We continue with the proof of Proposition 2.2. Let $\cJ \subset \cO_Z$
denote the (coherent) sheaf of ideals of the ramification locus $\br_h$.
For any  $k\in\N$ the  sheaf $\cS_k= \cJ^k \cdotp \cVT_Z$
(the product of $k$ copies of $\cJ$ with $\cVT_Z$) is
also coherent analytic. Fix $k$ and write $\cS=\cS_k$. Let $\Omega$ be a
Stein open subset of $Z$. We claim that there exist finitely many
sections $X_1,\ldots, X_N$ of $\cS$ over $\Omega$ which generate
$\cS$ at each point of $\Omega\bs \br_h$. This can be seen
by a standard argument as follows.

Since $\Omega$ is Stein, Cartan's Theorem A gives for each point
$z\in\Omega$ finitely many sections of $\cS$ over $\Omega$ which generate
$\cS$ at $z$. We can choose a point $a_j$ in each connected component 
$\Omega_j$ of $\Omega\bs \br_h$ such that the sequence $\{a_j\}$ is 
discrete in $\Omega$. Since $\cS$ is a free sheaf over the
complex manifold $\Omega\bs \br_h$, a simple (and well known)
extension of Cartan's Theorem A gives finitely many sections of $\cS$
over $\Omega$ which generate $\cS$ at each point $a_j$. The exceptional set
$Z_1$, consisting of points in $\Omega$ at which these
sections fail to generate $\cS$, is a complex subspace
of $\Omega$ satisfying $\dim Z_1\cap\Omega_j <\dim \Omega_j$
for each $j$. For each index $j$ for which $Z_1\cap \Omega_j \ne\emptyset$
we now choose a point $b_j\in Z_1\cap \Omega_j$ such that $\{b_j\}$ 
is discrete in $\Omega$. By Cartan's Theorem A there exist finitely 
many sections of $\cS$ over $\Omega$ which, 
together with the sections chosen in the first step, 
generate $\cS$ at each $b_j$. The exceptional set
$Z_2\subset \Omega$ at which all these sections fail to generate $\cS$
now satisfies $\dim Z_2\cap \Omega_j\le \dim \Omega_j -2$
for each $j$.

Continuing this way we obtain in finitely many steps sections
$X_1,\ldots, X_N$ of the sheaf $\cS$ over $\Omega$ which 
generate $\cS$ at each point of $\Omega\bs \br_h$.
(However, the minimal number of  generators of $\cS$ at ramification
points $z\in \br_h$ need not be bounded from above and hence
there need not exist finitely many sections generating $\cS$ over
$\Omega$.) By construction $X_1,\ldots, X_N$ are \hvf s
on $\Omega$ which are vertical with respect to $h$,
they vanish to order $k$ along $\br_h$, and they generate the vertical
tangent space $VT_z Z=\ker dh_z$ at each point $z\in \Omega\bs \br_h$.
Let $\phi^j_t$ denote the flow of $X_j$ for complex time $t$.
The map
$$
    s(z,t_1,\ldots,t_N)=\phi^1_{t_1}\circ \cdots \circ \phi^N_{t_N}(z),
$$
which is defined and \holo\ in an open \nbd\ of $\Omega\times \{0\}^N$
in $\Omega\times\C^N$ and takes values in $Z$,  satisfies Proposition 2.2.
Indeed, its partial derivative on $t_j$ at $t=0$ equals $X_j(z)$, and since
these vectors generate $VT_z(Z)$ for $z \in \Omega\bs \br_h$, $s$ satisfies
(iv). The properties (i)--(iii) are clear.
\endpr

\demo Proof of Theorem 2.1:
Proposition 2.2 enables us to prove Theorem 2.1 by following
step by step the proof of Theorem 1.4 in [FP3]. We shall point
out those places in the proof where a change or remark is needed.

The reader should first look at Theorem 5.2 in [FP3] (and Theorem 5.1 
in [FP1]). The situation is the following (we describe the basic 
case without parameters). We are given a Cartan pair $(A,B)$ in $X$, 
where the set $A$ contains a \nbd\ of the subvariety $X_0\subset X$
and where $h$ admits a dominating spray over a \nbd\ of $B$, 
and \holo\ sections $a$, $b$ of $h$ defined over a \nbd\ of $A$ resp.\ 
of $B$ such that $a$ and $b$ are uniformly close to each other 
over a \nbd\ of $C=A\cap B$. We wish to patch $a$ and $b$ into a 
single \holo\ section $\wt a$ over a \nbd\ of $A\cup B$ which is 
uniformly close to $a$ over $A$.  

The problem is reduced to the model situation
given by Proposition 5.2 in [FP1] (when $X_0=\emptyset$)
or by Proposition 4.2 in [FP3] (in the general case).
The model situation also applies in the present case without any
changes. The reduction is accomplished by Lemmas 5.3 and 5.4 
in [FP1]. Denote by $B^n(\e) \subset \C^n$ the open ball of 
radius $\e$ and center at the origin. These lemmas show
how to construct the following: 

\item{--} a local $h$-spray $s_1 \colon U_A \times B^n(\e)\to Z$ 
over a Stein \nbd\ $U_A \subset Z$ of $a(A)$ such that $s_1$
is dominating over a \nbd\ of $a(C)$,

\item{--} a global spray $s_2\colon U_B \times\C^n \to Z$
over a Stein \nbd\ $U_B\subset Z$ of $b(B)$ such that $s_2$ 
is dominating over a \nbd\ of $b(C)$, and

\item{--} an injective fiber preserving \holo\ map 
$\psi \colon \wt C\times B^n(\e)\to \wt C\times \C^n$
(where $\wt C\subset X$ is an open \nbd\ of $C=A\cap B$) such that
$$ s_2(b(x),\psi(x,t))= s_1(a(x),t) \qquad 
						(x\in \wt C,\ t\in B^n(\e)).
$$

The local spray $s_1$ is obtained from the spray 
$s$ granted by Proposition 2.2 above (the corresponding 
spray in [FP1] was denoted $\wt s$). The main point to 
observe with respect to the exposition in [FP1] is that 
the construction of $s_1$ only requires the domination property 
of $s$ over a \nbd\ of $a(C)$ (and not over a \nbd\ of $a(A)$).
Since the spray $s$ furnished by Proposition 2.2 is 
dominating outside of $\br_h$ and since 
$h(\br_h)\cap C=\emptyset$, $s$ is dominating over a \nbd\ 
of $a(C)$ as required. The spray $s_2$ is obtained 
from the global dominating spray over a \nbd\ of $B$ which exists 
by assumption. 

In order to get the transition map $\psi$ as above we must 
insure in addition that the kernels of $ds_1$ and $ds_2$ along the zero 
section are isomorphic (as \hvb s) over a \nbd\ $U_C\subset Z$
of $b(C)$ (which is chosen such that $U_C\subset U_A\cap U_B$).
The details of this construction are given by Lemmas 5.3 
and 5.4 in [FP1] (where $s_1$ and $s_2$ are constructed such 
that the above kernels are even close to each other and 
hence isomorphic). 

Since $X$ may have singularities contained in the subvariety
$X_0$, a remark is in order regarding the proof of Proposition 4.2
in [FP3] (the attaching lemma in the model situation). In the 
present case the patching is performed over open sets in 
$X\bs X_0 \subset X_{\rm reg}$. The $\dibar$-problems which 
arise in this patching have compact support contained in 
$X\bs X_0 \subset X_{\rm reg}$. Such $\dibar$-problems 
can be solved by transporting them to $\C^N$ via a 
holomorphic map $g\colon X \to \C^N$ which is a homeomorphism of $X$ 
onto a closed complex subvariety $\wt X\subset \C^N$ and which 
is biholomorphic on $X_{\rm reg}$. (Compare with section 7 in [FP3].)

We now proceed to section 6 of [FP3] where Theorem 1.4
of that paper is proved. The crucial step is furnished by Proposition 6.1
in [FP3]. To see that its proof remains valid in our current situation we 
observe that the sets $A_0,A_1,\ldots,A_n\subset X$,
which are chosen at the beginning of the proof of Proposition 6.1 in [FP3], are
such that $A_0$ contains a \nbd\ of $X_0$ while the sets $A_1,\ldots,A_n$ do not
intersect $X_0$. Since $h$ is a submersion of complex manifolds over $X\bs X_0$,
the techniques developed in [FP2, FP3] for holomorphic submersions onto Stein
manifolds can be applied whenever the first set $A_0$ is not involved.
Using those techniques we can patch any collection of holomorphic sections
$a_j\colon \wt A_j\to Z$ ($1\le j\le n$), where $\wt A_j \subset X\bs X_0$
is a small open \nbd\ of $A_j$ over which $h$ admits a spray, into a single
\holo\ section $b$ over a \nbd\  of $A^n=A_1\cup A_2\cup\ldots\cup A_n$, provided
that the sections $a_j$ belong to a holomorphic complex. (We are referring
to the transformation of a \holo\ complex associated to the Cartan string
$(A_1,\ldots,A_n)$ into a \holo\ section over their union $A^n$;
the details of this procedure are explained in [FP2, Proposition 5.1].)
The same procedure also gives a homotopy of \holo\ sections over a
\nbd\ of $A_0\cap A^n$ connecting $a$ and $b$.

It remain to patch $a$ and $b$ into a single \holo\ section over a
\nbd\ of $A_0\cup A^n$. This is accomplished as in [FP3] by combining the
homotopy version of the Oka-Weil approximation theorem (see e.g.\
Theorem 2.1 in [FP3]) with Theorem 5.2 in [FP3] which holds in the
current situation as explained above. The proof of Theorem 2.1
can now be concluded by the globalization procedure given
in [FP3] (proof of Theorem 1.4).

\beginsection 3.\ Proof of Theorem 1.3.

In this section we reduce Theorem 1.3 to Theorem 2.1.
Let $h\colon Z\to X$ and $f\colon Y \to X$ be as in Theorem 1.3. Set
$$
   \eqalignno{ \wh Z &= \{ (y,z)\colon y\in Y,\ z\in Z,\ f(y)=h(z)\}, \cr
   & \wh h(y,z)=y \in Y,\quad \sigma(y,z)=z\in Z. & (3.1) \cr}
$$
Clearly $\wh Z$ is a closed complex subspace of $Y\times Z$, the maps
$\wh h\colon \wh Z \to Y$ and $\sigma\colon\wh Z\to Z$ are \holo, and
we have $f\wh h=h\sigma$.

By assumption the set $Y_0=f^{-1}(X_0) \subset Y$ contains the singular 
locus $Y_{\rm sing}$. For each $y\in Y\bs Y_0$ we have $f(y)\in X\bs X_0$ 
and hence $h$ is a submersion over an open \nbd\ $U\subset X\bs X_0$ of $f(y)$. 
Setting $V=f^{-1}(U)$ it follows that $\wh h\colon \wh h^{-1}(V)\to V$ is a 
submersion. Thus $\wh h$ is a surjective submersion over $Y\bs Y_0$.
For any section $\wh g\colon Y\to\wh Z$ of $\wh h \colon \wh Z\to Y$
the map $g= \sigma \wh g\colon Y\to Z$ is a lifting of $f$
with respect to $h$:
$$
    hg=h(\sigma \wh g)=(h \sigma)\wh g
    =(f\wh h)\wh g = f(\wh h \wh g)= f.
$$
Moreover, any lifting $g$ of $f$ is of this form:
from $h(g(y))=f(y)$ ($y\in Y$) it follows that the point
$\wh g(y):=(y,g(y)) \in Y\times Z$ belongs to the
subset $\wh Z \subset Y\times Z$  (3.1) and
hence $\wh g\colon Y\to \wh Z$ is a section of $\wh h$. Furthermore,
$\sigma(\wh g(y))=\sigma(y,g(y))=g(y)$ whence $g$ is obtained from
the section $\wh g \colon Y\to \wh Z$. Therefore Theorem 1.3 follows
immediately from Theorem 2.1 and the following lemma.

%
%  Pulling back sprays
%
\proclaim 3.1 Lemma: {\bf (Pull-back sprays.)}
Let $f\colon Y\to X$ and $h\colon Z\to X$ be \holo\ maps.
Assume that $U\subset X$ is an open set such that
$h\colon Z|_U= h^{-1}(U) \to U$ is a submersion which
admits a spray. Then the map $\wh h\colon \wh Z\to Y$
defined by (3.1) is a submersion with spray over
$V=f^{-1}(U)\subset Y$.

\demo Proof:
Let $(E,p,s)$ be a spray associated to the submersion
$h\colon Z|_U\to U$ (Definition 1). Set $V=f^{-1}(U) \subset Y$
and observe that $\sigma$ maps $\wh Z|_V=\wh h^{-1}(V)$ to $Z|_U$.
Let $\wh p\colon \wh E\to \wh Z|_V$ denote the pull-back
of the \hvb\ $p\colon E\to Z|_U$ by the map
$\sigma \colon \wh Z|_V \to Z|_U$. Explicitly, we have
$$
   \eqalign{ \wh E &= \{ (\wh z,e)\colon
    \wh z\in \wh Z|_V,\ e\in E;\ \sigma(\wh z)= p(e)\} \cr
    &= \{(y,z,e)\colon y\in V,\ z\in Z,\ e\in E;\
       f(y)=h(z),\ p(e)=z\}; \cr
     & \!\!\!\! \wh p(\wh z,e) = \wh z.  \cr}
$$
Consider the map $\wh s\colon \wh E\to \wh Z|_V$, $s(y,z,e)=(y,s(e))$.
We claim that $(\wh E,\wh p,\wh s)$ is a spray associated to
the submersion $\wh h\colon \wh Z|_V\to V$.
We first check that  $\wh s$ is well defined.
If $(y,z,e)\in \wh E$ then $p(e)=z$ and $h(z)=f(y)$.
Since $s$ is a spray for $h$, we have $h(s(e))= h(z)=f(y)$
which shows that the point $\wh s(y,z,e)=(y,s(e)) \in Y\times Z$
belongs to the fiber $\wh Z_y$. This verifies property (i) in
Definition 1. Clearly $\wh s(y,z,0_{(y,z)})=(y,s(0_z))=(y,z)$ which
verifies property (ii) in Definition 1. It is also immediate
that $\wh s$ satisfies property (iii) provided that $s$
does since the vertical derivatives of the two maps
coincide under the identification $\wh Z_{y} \cong Z_{f(y)}$
and $\wh E_{(y,z)} \cong E_{z}$. This proves Lemma 3.1.

%
%
%  Multivalued sections...
%
%
\beginsection 4. Multivalued sections and analytic covers.

In this section we recall some well known results on symmetric products 
and multivalued sections which will be used in the proof of
Theorem 1.1. Our reference is Appendix V in [W].

Denote by $Z^d$ the $d$-fold Cartesian power of a set $Z$. 
The group $\Pi_d$ of all permutations on $d$ elements acts on 
$Z^d$ by permuting the entries, and we denote this action by $\rho$. 
The quotient space is called the {\bf $d$-fold symmetric power} of $Z$ 
and is denoted $Z^d_{\rm sym}$. For $z=(z_1,\ldots,z_d)\in Z^d$ we write 
$\pi(z)=[z]=[z_1,\ldots,z_d] \in Z^d_{\rm sym}$. 
A $d$-valued map from $X$ to $Z$ is a map $F\colon X\to Z^d_{\rm sym}$.
The number $d$ is called the {\bf degree} of $F$ and denoted $\deg F$. 

Assume from now on that $X$ and $Z$ are reduced complex spaces.
Then $\pi\colon Z^d\to Z^d_{\rm sym}$ induces a natural (quotient) 
complex structure on $Z^d_{\rm sym}$ such that \holo\ functions 
on $Z^d_{\rm sym}$ correspond to $\rho$-invariant \holo\ functions on 
$Z^d$. In particular, if $F=[f_1,\ldots,f_d] \colon X\to Z^d_{\rm sym}$
is a \holo\ map and if $P$ is a $\rho$-invariant holomorphic function on
$Z^d$ then $P(f_1,\ldots,f_d)$ is a well defined \holo\ function on $X$.

We recall some natural operations on symmetric products.

\item{1.} If $Z$ is a complex subspace of another complex space $S$
then $Z^d_{\rm sym}$ is in a natural way a subspace of 
$S^d_{\rm sym}$, and any map $X\to  S^d_{\rm sym}$ whose 
image belongs to $Z^d_{\rm sym}$ may also be considered 
as a map $X\to Z^d_{\rm sym}$. More generally, any \holo\ map 
$g\colon Z\to S$ induces a \holo\ map $\wh g\colon Z^d_{\rm sym}\to S^d_{\rm sym}$.

\item{2.} For any pair of integers $d,k\in \N$ we have a natural 
\holo\ map $\tau\colon Z^d_{\rm sym}\times Z^k_{\rm sym}\to Z^{d+k}_{\rm sym}$
induced by the identification $Z^d\times Z^k=Z^{d+k}$. 
Given a pair of maps $F_1\colon X\to Z^d_{\rm sym}$ and 
$F_2\colon X\to \times Z^k_{\rm sym}$, we denote 
$$ F_1\oplus F_2= \tau(F_1,F_2) \colon X\to Z^{d+k}_{\rm sym}. $$

The direct sum generalizes to several terms and we write 
$F=\oplus_j m_jF_j$, where the $F_j$'s are multivalued maps 
of $X$ to $Z$ and $m_j\in\N$. Clearly we have $\deg F=\sum_j m_j \deg F_j$. 
A map $F\colon X\to Z^d_{\rm sym}$ is called {\bf irreducible}
if it cannot be decomposed as a direct sum of multivalued maps 
of smaller degrees. We recall the notation 
$$ \mu_F= \max\{\#F(x)\colon x\in X\}, \quad
    \d_F=\{x\in X\colon \#F(x)<\mu_F\}, 
$$ 
where $\# F(x)$ denotes the number of distinct points 
in $F(x)$.  

\proclaim 4.1 Proposition: 
Assume that $F\colon X\to Z^d_{\rm sym}$ is a continuous 
(resp.\ holo) map such that $\d_F$ is nowhere dense in $X$ 
and $X\bs \d_F$ is pathwise connected and locally pathwise connected. 
Then $F$ has a decomposition $F=\oplus m_j F_j$ where the continuous 
(resp.\ \holo) maps $F_j\colon X\to Z_{\rm sym}^{d_j}$ are irreducible 
and the decomposition is unique up to the order of terms. 
Furthermore $\br F_j\subset \br F$ and $\d_{F_j}\subset \d_F$
for each $j$. Such a decomposition exists in particular if $X$ is 
an irreducible $n$-dimensional complex space and ${\cH}^{2n-1}(\d_F)=0$.

Proposition 4.1 is proved in [W, Appendix V]. The idea of the proof is 
as follows. Let $\wt V$ be the graph of $F$ over $\wt X=X\bs \d_F$,
i.e., $\wt V$ consists of all points in the fibers $F(x)$ for $x\in \wt X$.
Since $\wt X$ is connected and $\# F(x)=\mu_F$ for all $x\in \wt X$, 
$\wt V$ is a union of finitely many connected componets $\wt V_j$ such that
$h\colon \wt V_j\to \wt X$ is a finite unramified covering projection
onto $\wt X$, say with $d_j$ sheets, and $F$ has constant multiplicity 
$m_j\in \N$ along $\wt V_j$. Let $V_j$ denote the closure of $\wt V_j$ in $Z$. 
One can then show that $V_j$ is the graph of a $d_j$-valued section
$F_j$ of $h$ and $F=\oplus m_j F_j$.
\endpr

Assume now that $h\colon Z\to X$ is a surjective \holo\ map of reduced 
complex spaces. A map $F=[f_1,\ldots,f_d]\colon X\to Z^d_{\rm sym}$ 
such that $f_j(x)\in Z_x=h^{-1}(x)$ for all $x\in X$ and all $j$ is 
called a {\bf $d$-valued section} of $h$. The direct sum operation 
and Proposition 4.1 extend to multivalued sections. 

We shall assume from now on that $F\colon X\to Z^d_{\rm sym}$
is a $d$-valued section of $h\colon Z\to X$ which satisfies the
hypothesis of Proposition 4.1. If $F$ is irreducible, we define
its {\bf graph} $V(F)\subset Z$ by 
$$  
	V(F)=\{z\in Z \colon z\in F(x)\ {\rm for\ some}\ x\in X\}. 
$$
If $F=\oplus m_j F_j$ with $F_j$ irreducible for all $j$, 
we let $V(F)=\sum m_j V(F_j)$ be the disjoint union of $m_j$ 
copies of $V(F_j)$ for each $j$. If we consider $V(F)$ as a multiplicity 
subset of $Z$ then $h\colon V(F)\to X$ is a proper continuous map onto 
$X$ which is $d$-sheeted over $X\bs \d_F$.

Recall that an {\bf analytic chain} in $Z$ is a formal locally finite 
combination $V=\sum m_j V_j$ of closed complex subvarieties 
$V_j\subset Z$ with integer coefficients; if $m_j\ge 0$
for all $j$ then $V$ is called {\bf effective}. If $\dim V_j=n$ 
for all $j$ then $V$ is said to be purely $n$-dimensional. 
In the sense of currents we have $[V]=\sum_j m_j[V_j]$. 
The following proposition shows that there is a bijective correspondence
between multivalued \holo\ sections of $h\colon Z\to X$ and
analytic chains $V\subset Z$ such that $h|_V\colon V\to X$ is 
an analytic cover. Furthermore, when $X$ is irreducible, the decomposition 
of $F$ into irreducible components (given by Proposition 4.1) 
correponds to the decomposition of its graph $V(F)$ into irreducible 
complex subvarieties. 

\proclaim 4.2 Proposition:
Let $h\colon Z\to X$ be a \holo\ map of complex spaces.
Let $F\colon X\to Z^d_{\rm sym}$ be a \holo\ $d$-valued section of $h$ 
such that $X\bs \d_F$ is pathwise connected and locally pathwise connected
(this is always the case if $X$ is irreducible). If $F=\oplus m_j F_j$ 
is the decomposition into irreducible components granted by Proposition 4.1
then for each $j$ the graph $V(F_j)$ is a complex subvariety of $Z$ with 
finite proper $h$-projection onto $X$ (hence $V(F)=\sum m_j V(F_j)$ 
is an analytic chain in $Z$). Conversely, if $X$ is an irreducible 
$n$-dimensional complex space  and $V$ is an effective analytic chain in 
$Z$ such that $h|_V \colon  V\to X$ is a $d$-sheeted analytic cover 
then $V$ is the graph of a \holo\ $d$-valued section of $h$.

This can be proved by standard arguments from the theory of 
analytic covers. We omit the details and refer instead to [W]. 

The following lemma shows that local multivalued sections
of $h$ exist at each point where $h$ has maximal rank.

\proclaim 4.3 Lemma:
Let $h\colon Z\to X$ be a \holo\ map. Suppose that
$Z$ is locally irreducible at $z_0\in Z$, $X$ is locally
irreducible at $x_0:=h(z_0) \in X$ and
$$\dim_{z_0} h^{-1}(x_0)=\dim_{z_0}Z - \dim_{x_0}X.$$
Then there exist an integer $d\in\N$ and a local \holo\ $d$-valued
section $F$ of $h$ in a \nbd\ of $x_0$ such that $F(x_0)=[z_0,\ldots,z_0]$.

The number $k:=\dim_{z_0} h^{-1}(x_0)$ is called
the {\bf corank} of $h$ at $z_0$ and
$\dim_{z_0}Z-k$ is the {\bf rank} of $h$ at $z_0$.
Clearly the rank cannot exceed $\dim_{x_0}X$, and the
hypothesis in the lemma is that $h$ has maximal rank at $z_0$.

\demo Proof:
Since $h^{-1}(x_0)$ is a complex subvariety of $Z$ whose dimension at $z_0$
equals $k$, there exists a germ of an irreducible complex subvariety $V \subset Z$
at $z_0$ such that $\dim V+k=\dim_{z_0}Z$ and $z_0$ is an isolated point of
$V\cap h^{-1}(x_0)$. By a standard localization argument we obtain open \nbd s
$U\subset X$ of $x_0$ and $\wt U\subset Z$ of $z_0$ such that
$h\colon V\cap \wt U \to U$ is a proper finite map. The rank hypothesis
implies $\dim V=\dim_{x_0}X$ and hence $h(V\cap\wt U)=U$ provided
that $U$ is irreducible (as we may assume to be the case).
By Proposition 4.2 the set $V\cap \wt U$ is the graph of a
\holo\ multivalued section of $h$ over $U$.
\endpr

\beginsection 5. Proof of Theorem 1.1.

Let $F\colon X\to Z$ be a $d$-valued section of $h\colon Z\to X$
satisfying the hypotheses of Theorem 1.1. Thus $F$ is
\holo\ in an open set $U_0\subset X$ containing a
complex subvariety $X_0\subset X$ and unramified over $X\bs X_0$.
From $\cH^{2n-1}(\d_F)=0$ it follows by Proposition 4.1
that $F=\oplus m_j F_j$ for some irreducible multivalued sections
$F_j$ which are \holo\ over $U_0$, unramified over $X\bs X_0$
and satisfy $\cH^{2n-1}(\d_{F_j})=0$. It suffices to prove the result 
for each $F_j$. 

Thus we may assume without loss of generality that $F$ is an irreducible 
$d$-valued section satisfying $\br_F\subset X_0$, $\mu_F=d$ and
$\cH^{2n-1}(\d_F)=0$. The heart of the proof is the following lemma.

\medskip\ni\bf 5.1 Lemma. \sl
There exists a normal complex space $Y$ and a continuous map 
$g\colon Y \to Z$ such that, setting $f=hg\colon Y\to X$
and $Y_x=f^{-1}(x)$, we have:  
\item{(a)} $g$ is \holo\ in $f^{-1}(U_0)$ and $g(Y_x)=F(x)$ for each $x\in X$, 
\item{(b)} $f\colon Y\to X$ is a $d$-sheeted analytic cover which is 
unramified over $X\bs X_0$.
\medskip\rm

The proper way to think about $Y$ is as the `normalized graph' of $F$
where the self-intersections over $X\bs X_0$ have been removed.

\demo Proof: Over the set $U_0$ the $d$-valued section $F$ is \holo\ 
and hence its graph $V(F|_{U_0})$ is an effective chain in $h^{-1}(U_0)$ 
with finite proper $h$-projection onto $U_0$ (Proposition 4.2). Let 
$g_0 \colon Y_0\to V(F|_{U_0})$ denote its normalization (considered
as a map to $Z$). Then $g_0$ and $f_0\colon=hg_0\colon Y_0\to U_0$ 
satisfy the required properties over $U_0$. 

We now extend $(Y_0,g_0)$ as follows. Let $U\subset X\bs X_0$
be any open set such that $F|_U=\oplus_{j=1}^d F_j$ where $F_j\colon U\to Z$
are continuous sections of $h$ over $U$. From $\cH^{2n-1}(\d_F)=0$ 
it follows by Proposition 4.1 that the $F_j$'s are unique 
up to reordering (and there are no repetitions since $\mu_F=d$). Let 
$\N_d=\{1,2,\ldots,d\}$. Set $Y_U=U\times \N_d$ (the disjoint union
of $d$ copies of $U$) and define the map $g_U\colon Y_U \to Z$
by $g_U(x,j)=F_j(x)$. We introduce a complex structure on $Y_U$
by requiring that the (trivial) $d$-sheeted projection 
$f_U \colon= hg_U\colon Y_U\to U$ is biholomorphic on each sheet
$U\times\{j\}$.

If $U'\subset X\bs X_0$ is another open subset such that  
$F|_{U'}=\oplus_{j=1}^d F'_j$, it follows from $\cH^{2n-1}(\d_F)=0$ 
that for each connected component $\Omega$ of $U\cap U'$ there is a 
permutation $\sigma$ on $\N_p$ such that $F_j(x)=F'_{\sigma(j)}(x)$ 
for $x\in \Omega$ and $j=1,\ldots,d$. This defines a transition map 
$$ \sigma_{U,U'}\colon Y_U|_{\Omega} \to Y_{U'}|_{\Omega}, \quad 
            (x,j)\to (x,\sigma(j)).
$$
Clearly $\sigma_{U,U'}$ is biholomorphic with respect to the complex
structures on $Y_U$ and $Y_{U'}$ and $f_U=f_{U'}\circ \sigma_{U,U'}$.
Using these transition maps we may patch $Y_U$ and $Y_{U'}$ 
to a complex manifold $Y_{U\cup U'}$ which contains $Y_U$ and
$Y_{U'}$ as open subsets and such that the maps $g_U$ and $g_{U'}$
agree on the intersection of their domains to give a continuous map
$g_{U\cup U'}\colon Y_{U\cup U'}\to Z$. 

Since $F$ is unramified over $X\bs X_0$, we may globalize this
construction by covering $X\bs X_0$ with open sets as above and 
using the transition maps $\sigma_{U,U'}$ to construct a pair 
$(Y,g)$ with the required properties. 
\endpr

We continue with the proof of Theorem 1.1. Since $f\colon Y\to X$ is a 
finite map of $Y$ onto a Stein space $X$, it follows that the space $Y$ 
is also Stein [LeB]. The inverse $f^{-1}$ is a $d$-valued
\holo\ section of $f\colon Y\to X$. Consider the map $g\colon Y\to Z$ 
as a continuous lifting of $f=hg\colon Y\to X$ with respect to  
$h\colon Z\to X$. Theorem 1.3 furnishes a homotopy of liftings
$g_t\colon Y\to Z$ connecting $g_0=g$ to a holomorphic lifting $g_1$. 
Then $F_t=g_t f^{-1}$ ($0\le t\le 1$) is a homotopy of $d$-valued sections 
of $h$ satisfying Theorem 1.1.  
\endpr

\medskip
\ni \bf A remark on [FP3]. \rm We take this occasion to point
out the following omission in the hypothesis of Theorem 5.2
in [FP3]: {\it The sets $A$ and $B$ in the statement of that theorem
must have a basis of open Stein \nbd s in $X$}. (This property was
assumed in the closely related Theorem 5.1 in [FP3], but was accidentally
omitted in Theorem 5.2.) The reader can observe that in all applications
of Theorem 5.2 in [FP3] this additional hypothesis is satisfied.

\medskip
\ni\bf Acknowledgements. \rm  I wish to thank T.\ Ohsawa
who raised the question answered in part by Theorems 1.1 and 1.3
at the Hayama Symposium on Complex Analysis 2000. 
I wish to thank J.\ Prezelj for her interest and stimulating 
discussions on this subject. This research was supported in part 
by an NSF grant and a grant from the Ministry of Science of the 
Republic of Slovenia.

%
%
%  References
%
%
\medskip\ni\bf References. \rm

\ii{[BV]} W.\ Barth, P.\ C.\ Van de Ven:
Compact complex surfaces. 
Ergebnisse der Mathematik und ihrer Grenzgebiete (3), 4. 
Springer-Verlag, Berlin, 1984.

\ii{[Fi]} G.\ Fischer: Complex analytic geometry. 
Lecture Notes in Mathematics, Vol.\ 538. Springer,
Berlin-New York, 1976.

\ii{[FR]} O.\ Forster and K.\ J.\ Ramspott:
Analytische Modulgarben und Endromisb\"undel.
Invent.\ Math.\ {\bf 2}, 145--170 (1966).

\ii{[F1]} F.\ Forstneri\v c: On complete intersections. 
Ann. Inst. Fourier {\bf 51} (2001), 497--512.

\ii{[F2]} F.\ Forstneri\v c: 
The Oka principle, lifting of holomorphic maps and removability of intersections. 
Proc.\ of Hayama Symposium on Several Complex Variables 2000,
pp.\ 49--59, Japan, 2001.

\ii{[F3]} F.\ Forstneri\v c: The Oka principle for sections of
subelliptic submersions. Preprint, 2001.

\ii{[FP1]} F.\ Forstneri\v c and J.\ Prezelj:
Oka's principle for holomorphic fiber bundles with sprays.
Math.\ Ann.\ {\bf 317} (2000), 117-154.

\ii{[FP2]} F.\ Forstneri\v c and J.\ Prezelj:
Oka's principle for holomorphic submersions with sprays.
Math.\ Ann., to appear.

\ii{[FP3]} F.\ Forstneri\v c and J.\ Prezelj:
Extending holomorphic sections from complex subvarieties.
Math.\ Z.\ {\bf 236} (2001), 43--68.

\ii{[Gr1]} H.\ Grauert:
Holomorphe Funktionen mit Werten in komplexen Lieschen Gruppen.
Math.\ Ann.\ {\bf 133}, 450--472 (1957).

\ii{[Gr2]} H.\ Grauert: Analytische Faserungen \"uber
holomorph-vollst\"andigen R\"aumen.
Math.\ Ann.\ {\bf 135}, 263--273 (1958).

\ii{[Gro]} M.\ Gromov:
Oka's principle for holomorphic sections of elliptic bundles.
J.\ Amer.\ Math.\ Soc.\ {\bf 2}, 851-897 (1989).

\ii{[GR]} C.\ Gunning, H.\ Rossi:
Analytic functions of several complex variables.
Prentice--Hall, Englewood Cliffs, 1965.

\ii{[L]} F.\ L\'arusson: 
Excision for simplicial sheaves on the Stein site
and Gromov's Oka principle. Preprint, December 2000.

\ii{[LeB]} P.\ Le Barz: A propos des rev\^etements d'espaces de Stein.
Math.\ Ann.\ {\bf 222}, 63--69 (1976).

\ii{[O]} K.\ Oka: Sur les fonctions des plusieurs variables. III:
Deuxi\`eme probl\`eme de Cousin.
J.\ Sc.\ Hiroshima Univ.\ {\bf 9} (1939), 7--19.

\ii{[S]} Y.T.\ Siu: Every Stein subvariety admits a Stein neighborhood.
Invent.\ Math.\ {\bf 38}, 89--100 (1976).

\ii{[Wd]} G.\ W.\ Whitehead: Elements of homotopy theory.
Graduate Texts in Mathematics, 61. Springer, New York--Berlin, 1978.

\ii{[W]} H.\ Whitney: Complex Analytic Varieties.
Addison-Wesley, Reading, 1972.

\bigskip\medskip
\ni\bf Address: \rm Institute of Mathematics, Physics and Mechanics,
University of Ljubljana, Jadranska 19, 1000 Ljubljana, Slovenia

\bye